\documentclass[12pt,a4paper]{amsart}

\usepackage{amsmath}
\usepackage{amsthm,amssymb,amsfonts}

\usepackage{graphicx}

\newtheorem{thm}{Theorem}[section]
\newtheorem{lemma}[thm]{Lemma}
\newtheorem{prop}[thm]{Proposition}

\newtheorem{example}[thm]{Example}
\newtheorem{ques}[thm]{Question}

\theoremstyle{definition}
\newtheorem{example2}[thm]{Example}

\def\eqn#1$$#2$${\begin{equation}\label#1#2\end{equation}}

\def\th{{\leavevmode\setbox1=\hbox{t}
  \hbox to \wd1{t\kern-0.6ex{\char039}\hss}}}
\def\dh{{\leavevmode\setbox1=\hbox{d}
  \hbox to 1.05\wd1{d\kern-0.4ex{\char039}\hss}}}
\def\=#1{\if #1u{\accent23u}\else
\ifx #1d{\dh}\else \ifx #1t{\th}\else
 {\accent20 #1}\fi\fi\fi}
\def\'#1{\if #1i{\accent19\i}\else {\accent19 #1}\fi}

\def\ce{\mathbb C}

\def\en{\mathbb N}
\def\er{\mathbb R}

\let\Re\re

\def\dist{\operatorname{dist}}

\def \Ker {\operatorname{Ker}}

\def \Int {\operatorname{Int}}

\def\sp{\operatorname{span}}

\def \reg {\partial _{\kern1pt\text{reg}}}

\def\wscl#1{\overline{#1}^{w^*}}
\def\ncl#1{\overline{#1}^{\|\cdot\|}}

\def\Dd#1#2{{#2}^{(#1)}}
\def\Dt#1#2{{#2}^{(\widetilde{#1})}}
\newtoks\by
\newtoks\paper
\newtoks\book
\newtoks\jour
\newtoks\yr
\newtoks\pages
\newtoks\vol
\newtoks\publ
\newtoks\eds
\newtoks\proc
\newtoks\mathrev
\newtoks\web
\def\ota{{\hbox{???}}}
\def\cLear{\by=\ota\paper=\ota\book=\ota\jour=\ota\yr=\ota
\pages=\ota\vol=\ota\publ=\ota}
\def\endpaper{\the\by, \textit{\the\paper},
{\the\jour} \textbf{\the\vol} (\the\yr), \the\pages.\cLear}
\def\endbook{\the\by, \textit{\the\book}, \the\publ, \the\yr.\cLear}
\def\endprep{\the\by, \textit{\the\paper}, \the\jour.}
\def\endprepkma{\the\by, \textit{\the\paper}, \the\jour, available on http://adela.\-karlin.\-mff.\-cuni.\-cz/\~{ }rokyta/\-preprint/\-index.php.\cLear}
\def\endproc{\the\by, \textit{\the\paper}, \the\book,
\the\publ, \the\yr, \the\pages.\cLear}
\def\endper{\the\by, \textit{personal communication}.\cLear}

\newcommand{\na}{\mathbb{N}}
\newcommand{\cn}{\mathbb{C}}
\newcommand{\cu}{\mathbb{K}}

\newcommand{\ic}[1]{\operatorname{Int}_{\|.\|}\wscl{#1}}
\newcommand{\icu}{\ic{U}}

\title[Envelopes of open sets]
{Envelopes of open sets and extending holomorphic functions on dual Banach spaces}

\author[D. Garc\'{\i}a]{Domingo Garc\'{\i}a }
    \address{Departamento de An\'{a}lisis \\
Matem\'{a}tico\\ Universidad de Valencia\\  Valencia, Spain}
 \email{Domingo.Garcia@uv.es}

\author[O.F.K. Kalenda]{Ond\v{r}ej  F.K. Kalenda }
    \address{Department of Mathematical Analysis \\
Faculty of Mathematics and Physic\\ Charles University\\
Sokolovsk\'{a} 83, 186 \ 75\\Praha 8, Czech Republic}
 \email{kalenda@karlin.mff.cuni.cz}

\author[M. Maestre]{Manuel Maestre}
    \address{Departamento de An\'{a}lisis \\
Matem\'{a}tico\\ Universidad de Valencia\\  Valencia, Spain}
 \email{Manuel.Maestre@uv.es}

\thanks{The first and the third authors were  supported in part
by  MEC and FEDER Project MTM2008-03211. The third author was also supported by Prometeo 2008/101}
\thanks{The second author was supported in part by the grant
GAAV IAA 100190901  and in part by the Research Project
MSM~0021620839 from the Czech Ministry of Education.}

\date{}


\begin{document}

\begin{abstract}
We investigate certain envelopes of open sets in dual Banach spaces which are related to extending
holomorphic functions. We give a variety of examples of absolutely convex sets showing that the
extension is in many cases not possible. We also establish connections to the study of iterated
weak* sequential closures of convex sets in the dual of separable spaces.
\end{abstract}

    \maketitle

\section{Introduction and Preliminaries}

In \cite{CaGaMa}, Carando et al. studied theorems of Banach-Stone type for some algebras of
holomorphic functions on Banach spaces. The case of polynomials was studied, for instance, by
Cabello S\'{a}nchez, Castillo and R. Garc\'\i a  \cite{CCG},  D\'{\i}az and Dineen \cite{DD} and Lassalle
and  Zalduendo  \cite{LZ}. The authors in \cite{CaGaMa} proved, for example, that whenever you
have two Banach spaces $X$ and $Y$, one of them having the approximation property, two absolutely
convex  open sets $U\subset X^\ast$ and $V\subset Y^\ast$ and the algebras $H_{w^\ast u}(U)$ and
$H_{w^\ast u}(V)$, of holomorphic functions which are uniformly $w(X^\ast,X)$-continuous on
$U$-bounded (resp. $V$-bounded) sets, see below for the definitions,  are topologically algebra
isomorphic, then $X$ and $Y$ are isomorphic. In this study a key element is the description of the
spectra of that algebras. Some questions remained open and we are going to deal with them. Before
stating these questions we need to introduce some notations.

We refer to \cite{D} for background information on infinite dimensional complex analysis. In what follows $X$ and $Y$ will stand for real or complex Banach spaces. The motivation for our investigation comes from complex spaces but many results are true also for the real ones (with a simpler proof, usually). So, $\cu$ will denote the respective field -- either $\ce$ or $\er$. 
 By $B_X$ we will denote the open unit ball of $X$, $D_X$ will denote the closed unit ball of $X$.

Let $U\subset X$ be an open set.
 A subset
$B$ of $U$ is \emph{$U$-bounded} if it is bounded and has positive distance  to $X\setminus U$. A
countable family $\mathcal{B}=(B_n)_{n=1}^\infty$ is  a \emph{fundamental sequence of $U$-bounded sets}
if it satisfies the following conditions:
\begin{itemize}
 \item[(i)] $B_n$ is $U$-bounded for all $n$;
 \item[(ii)] for any $U$-bounded set $B$ there
exists $n$ such that $B\subset B_n$;
\item[(iii)] there exists a sequence of positive numbers
$(r_n)_{n=1}^\infty$ such that $B_n+r_n B_X\subset B_{n+1}$ for all $n$.
\end{itemize}
It is easy to check that each open set $U$ admits a fundamental sequence of $U$-bounded sets.  The canonical choice is $(U_n)_{n=1}^\infty$ where
$$U_n=\left\{x \in U : \|x \|\leq n \mbox{ and
} {\rm dist}(x ,X\setminus U)\geq\frac{1}{n}\right\},\quad n\in \mathbb{N}.$$
This family will be addressed as \emph{the canonical fundamental family of $U$-bounded sets}.

For an open subset $U$ of a complex space $X$, $H_{wu}(U)$ will denote the space of holomorphic functions on $U$
that are uniformly weakly continuous on $U$-bounded sets. It is  endowed with the topology of the
uniform convergence on  $U$-bounded sets, i.e. with the family of seminorms
$$\|f\|_B:=\sup\{|f(x)| : x\in B\}, \qquad B\mbox{ $U$-bounded}.$$
This space is a Fr\'{e}chet algebra since it is a closed subalgebra of the space of holomorphic
functions on $U$, bounded on $U$-bounded sets, which is Fr\'{e}chet when endowed with that topology
(see e.g. \cite{Barroso}).

 The spectrum of this algebra, consisting in all continuous
homomorphisms from $H_{wu}(U)$ into $\mathbb{C}$, will be denoted by
$M_{wu}(U)$.

Analogously, for an open subset $U$ of $X^\ast$, $H_{w^\ast u}(U)$ will stand the space of
holomorphic functions on $U$ which are uniformly $w(X^\ast,X)$-continuous on $U$-bounded sets. As
in the case of $H_{wu}(U)$, this space is also a Fr\'{e}chet algebra endowed with the topology of the
uniform convergence on $U$-bounded sets. The spectrum of this algebra will be denoted by
$M_{w^\ast u}(U)$.

The most obvious homomorphism in the spectra is the evaluation at a point $x\in U$. When $U=B_X$,
it is well known  that by uniform continuity any function $f\in H_{wu}(B_X)$ has a unique
extension $\tilde{f}$ to $B_{X^{\ast \ast}}$ that is uniformly weak-star continuous on $rB_{X^{\ast \ast}}$
for all $0<r<1$. But, by \cite[Theorem 5]{DaGa}  the function $f$ has an holomorphic  extension to $B_{X^{\ast
\ast}}$. Hence, by the uniqueness  of the continuous extension
 $\tilde{f}\in H_{w^\ast u}(B_{X^{\ast \ast}})$. So to each point $x^{\ast
\ast}\in B_{X^{\ast \ast}}$ we can associate a continuous homomorphism $\delta_{x^{\ast \ast}}\in
M_{wu}(B_X)$ by $\delta_{x^{\ast \ast}}(f)=\tilde{f}(x^{\ast \ast})$ for every $f\in H_{wu}(B_X)$.
The same happens when $U$ is an absolutely convex  open set, replacing $B_{X^{\ast \ast}}$ by the
$\|\cdot\|$-interior in $X^{\ast \ast}$ of the $w(X^{\ast \ast},X^{\ast})$-closure of $U$ (see
\cite[Remark 13]{CaGaMa}).

Thus to move around this stuff, finding a description of the
spectra, is natural to start in the dual $X^\ast$ of $X$.

Let $U\subset X^\ast$ be an open set and $(U_n)_{n=1}^\infty$ a fundamental family of $U$-bounded sets. Then,
every $f\in H_{w^\ast u}(U)$ is uniformly weak-star continuous in
each $U_n$, so it extends uniquely to a weak-star continuous
functions $\tilde{f}$ on $\wscl{U_n}$ with $\|\tilde f\|_{\wscl{U_n}}=\| f\|_{U_n}$. Hence each $x^\ast\in\wscl{U_n}$ produces a continuous homomorphism on $H_{w^\ast u}(U)$.

We set
$$\tilde U:=\bigcup_n \wscl{U_n}.$$
This set
coincides with
$$\bigcup \{\wscl B : B \text{ is } U\text{-bounded}\}$$
and  it is an open subset of $X^\ast$ contained in $\icu$, the norm-interior of
$\overline{U}^{w^\ast}$. A natural mapping can be defined:
$$\delta:\tilde U\longrightarrow M_{w^\ast u}(U), \ \ \ \delta(x^\ast):=\delta_{x^\ast}.$$

In \cite{CaGaMa} some characteristics of the relationship between
$U$, $\tilde U$ and $M_{w^\ast u}(U)$ were pointed out. The main questions concerning the relationship between $U$ and $\tilde U$ which we address also in this paper are the following:
\begin{enumerate}
    \item When $\tilde U=\icu$?
    \item Let $(U_n)_{n=1}^\infty$ be a fundamental family of $U$-bounded sets. Is $(\wscl{U_n})_{n=1}^\infty$ a fundamental family of $\tilde U$-bounded sets?
\end{enumerate}

The sets $U$ for which the second question has positive answer are called \emph{boundedly regular} in \cite{CaGaMa}.
It was noticed in \cite[Remark 5(ii)]{CaGaMa} that the answer to both questions is positive if $U$ is bounded and absolutely convex.

The importance of the second question is clear from \cite[Proposition 3]{CaGaMa} which we now recall.
Trying to clarify the properties of the extension of any element of $H_{w^\ast
u}(U)$, in \cite{CaGaMa} was introduced a new class of Fr\'{e}chet
algebras. For an open subset $U$ of $X^\ast$, let
$\mathcal{B}=(B_n)_{n=1}^\infty$ be a countable family of $U$-bounded sets
whose union is $U$  and satisfying that for each $n$ there is $r_n>0$ with $B_{n}+r_n B_{X^\ast}\subset B_{n+1}$. (Note that this is a weakening of the notion of fundamental family of $U$-bounded sets.) Endowed with the topology of the uniform convergence on these sets,
$$H_{\mathcal{B}w^\ast u}(U):=\{f\in H(U) : f|_{B_n} \text{ is weak-star uniformly continuous on } B_n \text{ for all } n\in \mathbb{N}\}$$
is a Frech\'{e}t algebra. In the case that $\mathcal{B}$ be a fundamental sequence of
$U$-bounded sets it follows that the equality $H_{\mathcal{B}w^\ast
u}(U)= H_{w^\ast u}(U)$ holds algebraically and topologically.
Now we are ready to recall the promised proposition from \cite{CaGaMa} which is a generalization of \cite[Theorem 3]{AR}.

\begin{prop}  \label{espectrow} {\rm (\cite[Proposition 3]{CaGaMa})} Let $X$ be a complex Banach space.
\begin{itemize}
\item[(i)] Let $U$ be a  balanced open subset of $X^\ast$, $(U_n)_{n=1}^\infty$ a fundamental family of $U$-bounded sets and
$\mathcal{B}=(\overline {U_n}^{w^\ast})_{n=1}^\infty$. Every $f\in
H_{w^\ast u}(U)$ extends uniquely to an $\tilde{f} \in
H_{\mathcal{B}w^\ast u}(\tilde U)$ and the mapping $i:H_{w^\ast
u}(U)\longrightarrow H_{\mathcal{B}w^\ast u}(\tilde U)$,
$i(f):=\tilde f$  is a topological algebra isomorphism. \noindent

\item[(ii)] If $U$ is an absolutely convex  open subset of
$X^\ast$ and $X$ has the approximation property then $M_{w^\ast u
}(U)=\delta(\tilde U)$.
\end{itemize}
\end{prop}

If $U$ is boundedly regular, it is clear that in the assertion (i) of this proposition we can substitute $H_{\mathcal{B}w^\ast u}(\tilde U)$ by $H_{w^\ast u}(\tilde U)$, hence the algebras $H_{w^\ast u}(U)$ and $H_{w^\ast u}(\tilde U)$ are canonically isomorphic.

In \cite[Example 7]{CaGaMa} an unbounded balanced open set $U$ was constructed such that $U$ is not boundedly regular and, moreover, the respective Fr\'echet algebras $H_{w^\ast u}(U)$ and $H_{w^\ast u}(\tilde U)$ are different.
In the present paper we complete this example by  some more counterexamples to the above questions.

In Section~\ref{s-balanced} we deal with bounded balanced sets. We prove that there is a bounded balanced open set $U$ such that $\tilde U=U\ne \icu$.
Note that this $U$ is boundedly regular for trivial reasons.

In Section~\ref{s-convex} we study the envelopes of open convex sets. It turns out that the convex case is very different from the balanced case. In particular, if $U$ is convex and boundedly regular, then necessarily $\tilde U=\icu$.

In Section~\ref{s-examples} we use a method developped in the previous section to construct several counterexamples. We show that there are unbounded absolutely convex open sets $U$ such that $\tilde{\tilde U}\ne \tilde U$. In fact, the iteration of the $\widetilde{\quad}$ operation may produce long transfinite one-to-one sequences. We also prove that an absolutely convex open set $U$
 need not be boundedly regular even if $\tilde{\tilde U}=\tilde U$.

In Section~\ref{s-algebras}  we prove that for absolutely convex $U$ the Fr\'echet algebras $H_{w^\ast u}(U)$ and $H_{w^\ast u}(\tilde U)$ are canonically isomorphic if and only if $U$ is boundedly regular.
 Thus for the construced examples the respective Fr\'echet algebras are different.

In the last section we collect some open questions.

\section{Examples of bounded balanced sets}\label{s-balanced}

The aim of this section is to prove the following example.

\begin{example} There is a separable Banach space $X$ and a balanced open set $U\subsetneqq B_{X^*}$ such that $\tilde U=U$ and $U$ is weak* dense in $B_{X^*}$. We can take, for example, $X=\ell_1$.
\end{example}

The following theorem specifies for which spaces $X$ we are able to find such an example in $X^*$. The previous example follows immediately from Theorem~\ref{bal-cpcp} and Example~\ref{exaellinfty}.

\begin{thm}\label{bal-cpcp} Let $X$ be a separable Banach space. Suppose there is $\varepsilon>0$ such that each nonempty relatively weak* open subset of $D_{X^*}$ has norm-diameter at least $\varepsilon$. Then there is a balanced open set $U\subsetneqq B_{X^*}$ such that $\tilde U=U$ and $U$ is weak* dense in $B_{X^*}$.
\end{thm}

Before proving the theorem let us discuss which spaces satisfy the assumptions.
First note that $X^*$ must be nonseparable. Indeed, if $X^*$ is separable, then any nonempty bounded subset of $X^*$ admits nonempty relatively weak* open subsets (in fact slices) of arbitrarily small diameter. This follows for example from \cite[Theorem 2.19 and Lemma 2.18]{phelps}. An easy positive example is the following one.

\begin{example}\label{exaellinfty} Any nonempty weak* open subset of the unit ball of $\ell_\infty=\ell_1^*$ has diameter $2$.
\end{example}

\begin{proof} Note that the closed unit ball of $\ell_\infty$ is
$$D=\{t:|t|\le1\}^\en$$
and the weak* topology coincides with the product topology. Let $U$ be a nonempty open subset of
$D$. Choose $(x_n)_n\in U$. Then there is $N\in \en$ and $\varepsilon>0$ such that
$$\{ (y_n)_n\in D: |y_n-x_n|<\varepsilon \mbox{ for }n=1,\dots,N\}
\subset U.$$
Define two points, $y=(y_n)_n$ and $z=(z_n)_n$ as follows:
\begin{gather*}
 y_n=z_n=x_n,\qquad n=1,\dots,N, \\ y_n=1, z_n=-1,\qquad n>N.
\end{gather*}

Then $y,z\in U$ and $\|y-z\|=2$.
\end{proof}

We may ask whether there are some other spaces satisfying the assumptions of
Theorem~\ref{bal-cpcp}. There are some results in this direction. For example,
in \cite[Theorem 2.5]{BLR} it is proved that whenever $A$ is an infinite-dimensional $C^*$-algebra, then every nonempty relatively weakly open subset of the unit ball has diameter $2$. If $A$ is moreover a dual space, then a fortiori nonempty relatively weak* open subsets of the unit ball have diameter $2$. So, this covers our example $\ell_\infty=\ell_1^*$ and also the noncommutative version $B(\ell_2)=N(\ell_2)^*$, the space of bounded linear operators on $\ell_2$ which is the dual to the space of all nuclear operators on $\ell_2$. Some generalizations of the results of \cite{BLR} are contained in the recent paper \cite{AcosBe}.

There are also some related results
in the realm of real spaces -- see \cite[Section III.1]{DGZ}. We are convinced that at least some
of them can be proved for complex spaces as well, but we have not checked it.

Now we give a series of lemmas which prove Theorem~\ref{bal-cpcp}.
We start by the following lemma on ``cones'' which we will need to compute $\tilde U$ for certain sets $U$.

\begin{lemma}\label{cones} Let $X$ be a  Banach space. Let $x\in X\setminus \{0\}$ and $0<a<\|x\|$.
Set
$$ C:=\bigcup\{t(x+a B_X), \ \ t\in \mathbb{K}\setminus \{0\}\}$$ $$
D:=\bigcup\{t(x+a D_X), \ \ t\in \mathbb{K}\}. $$ Then $C=\Int_{\|\cdot\|}D$.

\end{lemma}

\begin{proof}
Without loss of generality we may assume $\|x\|=1$. Let $y\in X$ arbitrary, $y\notin \sp \{x\}$.
Set $E:=\sp\{x,y\}$. Let $\varphi\in E^*$ such that $\varphi (x)=1$, $\varphi(y)=0$. Then
\begin{itemize}
    \item [(1)] $y\in C$ if and only if $\|\varphi\|>\frac{1}{a}$
    \item [(2)] $y\in D$ if and only if $\|\varphi\|\geq \frac{1}{a}$ \ .
\end{itemize}
Indeed, if $y\in C$ then there is $t\in \mathbb{K}\setminus \{0\}$ such that $\|y-tx\|<|t|a$.
Hence $|\varphi(y-tx)|=|t|>\frac{1}{a}\|y-tx\|$. Thus $\|\varphi\|>\frac{1}{a}$. Analogously, if
$y\in D$ then there is $t\in \mathbb{K}$ such that $\|y-tx\|\leq |t|a$. So $|\varphi(y-tx)|\geq
\frac{1}{a}\|y-tx\|$. As $y-tx\neq 0$ ($y\notin\sp\{x\}$), we get $\|\varphi\|\ge\frac1a$.

If  $y\notin C$ then   $\|y-tx\|\geq |t|a$ for all $t\in \mathbb{K}\setminus \{0\}$. For $t=0$
this inequality also holds. Hence $\|\alpha y-tx\|=|\alpha| \|y-\frac{t}{\alpha}x\|\geq |t|a$
for $\alpha\neq 0$. This also remains true for $\alpha =0$. Therefore $|\varphi (\alpha
y-tx)|=|t|\leq \frac{1}{a}\|\alpha y-tx\|$. Thus $\|\varphi\|\leq \frac{1}{a}$. Analogously,
if $y\notin D$ then
 $\|y-tx\|> |t|a$ for all $t\in \mathbb{K}$.  Hence   $|\varphi (\alpha y-tx)|< \frac{1}{a}\|\alpha
y-tx\|$ for $\alpha y-tx\neq 0$. Thus $\|\varphi\|< \frac{1}{a}$ (note that $\varphi$ is defined on a two-dimensional space and hence attains its norm).

Now, we  prove that $C=\Int_{\|\cdot\|}D$. As $C$ is open and $C\subset D$, we have $C\subset\Int_{\|\cdot\|} D$. We will show the converse inclusion.

Suppose $y\notin C\cup \{0\}$. Then $y\notin\sp\{x\}$ and we can define $E$ and $\varphi$ as above. Moreover, we have that
$\|\varphi\|\leq \frac{1}{a}$. Let $\psi\in E^*$ such that $\psi (x)=1$, $\|\psi\|=1$. Set
$\eta_\delta = (1-\delta)\varphi + \delta \psi$, $\delta\in (0,1)$. Then $\eta_\delta(x)=1$
and $\|\eta_\delta\|\leq (1-\delta)\frac{1}{a}+\delta<\frac{1}{a}$ (recall that $a<1$). Let $z\in\Ker \psi$, $z\neq 0$. If $\varphi(z)=0$, then $\Ker\psi=\Ker\varphi$ and hence $\psi$ is a multiple of $\varphi$. As $\varphi(x)=\psi(x)=1$, we conclude $\varphi=\psi$, so $\|\varphi\|=1<\frac1a$. Thus $y\notin D$.

Now suppose that $\varphi(z)\neq 0$. Set
$$\mu_\delta:=y-\frac{\delta}{1-\delta}\frac{\psi(y)}{\varphi(z)}z \ .$$
Then
$\mu_\delta\rightarrow y$ for $\delta\rightarrow 0^+$ and
$$\eta_\delta(\mu_\delta)=(1-\delta)\varphi(\mu_\delta) + \delta
\psi(\mu_\delta)=(1-\delta)\left(-\frac{\delta}{1-\delta}\frac{\psi(y)}{\varphi(z)}\varphi(z)\right)+\delta\psi(y)=0.$$
Since $\|\eta_\delta\|<\frac{1}{a}$ it follows that $\mu_\delta\notin D$.
Thus $y\notin \Int_{\|\cdot\|}D$.

Finally, suppose that $y=0$. Let $\psi\in X^*$ such that  $\|\psi\|=\psi(x)=1$ and  choose
$z\in\Ker\varphi$, $z\ne 0$. As $\|\psi\|=1<\frac{1}{a}$, we have $z\notin D$. Hence $tz\notin D$
for all $t\in (0,1)$ and $tz\rightarrow 0=y$ for $t\rightarrow 0^+$. So, again $y\notin
\Int_{\|\cdot\|}D$ and the proof is completed.
\end{proof}

\begin{lemma}\label{cutcones} Let $X$, $x$, $a$, $C$ and $D$ have the same meaning as in Lemma~\ref{cones}. Set
$$C'=C\cap B_{X},\qquad D'=D\cap D_{X}.$$
Then
\begin{itemize}
    \item[(a)] $D'\subset \bigcup \{t(x+a D_X), \ \ t\in \mathbb{K},|t|\le \frac{1}{\|x\|-a}\}$.
    \item[(b)] $\Int_{\|\cdot\|} D'=C'$.
    \item[(c)] $\ncl{C'}=D'$. If $X$ is, moreover, a dual Banach space equipped with the dual norm, then $D'$ is weak* closed and hence $\wscl{C'}=D'$.
\end{itemize}
\end{lemma}

\begin{proof}
(a) Let $y\in D'$. As $y\in D$, there is $t\in\cu$ such that $y=t(x+a D_X)$. Hence there is $z\in D_X$ with $y=t(x+a z)$. If $|t|>\frac{1}{\|x\|-a}$, then
$$\|y\|\ge|t|(\|x\|-a\|z\|)\ge|t|(\|x\|-a)>1,$$
so $y\notin D'$. This proves (a).

(b) This follows from Lemma~\ref{cones}, as
$$\Int_{\|\cdot\|} D'=\Int_{\|\cdot\|} (D\cap D_{X})=\Int_{\|\cdot\|} D\cap \Int_{\|\cdot\|} D_X=C\cap B_X=C'.$$

(c) Let $y\in D'$. By (a) there is some $t\in\cu$ with $|t|\le\frac{1}{\|x\|-a}$ and $z\in D_X$ with $y=t(x+a z)$. For $k\in\en$ set
$$y_k=\frac{t}{1+\frac{2a}{k(\|x\|-a)}}\cdot\left(x+a\left(1-\frac1k\right)z\right).$$
It is clear that $y_k$ norm-converges to $y$. Moreover, $y_k\in C'$. Indeed, obviously $y_k\in C$
and, moreover,
$$\|y_k\|\le \frac{\|y\|+\|t a \frac z k\|}{1+\frac{2a}{k(\|x\|-a)}}
\le \frac{1+\frac{a}{k(\|x\|-a)}}{1+\frac{2a}{k(\|x\|-a)}}<1.$$
This proves that $C'$ is norm-dense in $D'$.

It remains to prove that $D'$ is closed in the norm topology (in the weak* topology if $X$ is a dual space). We will prove it simultaneously:

Let $y$ be in the closure of $D'$. Then there is a net $y_\tau$ in $D'$ converging to $y$. For each $\tau$ there is $t_\tau\in \cu$, $|t_\tau|\le \frac{1}{\|x\|-a}$ and $z_\tau\in D_X$ with $y_\tau=t_\tau(x+a z_\tau)$.
As the net $t_\tau$ is bounded, we can assume (up to passing to a subnet) that it converges (in $\cu$) to some $t$. If $t=0$, then $y_\tau\to 0$, hence $y=0\in D'$. So suppose that $t\ne 0$. Then we can suppose that $t_\tau\ne0$ for each $\tau$. Therefore,
$$z_\tau=\frac{y_\tau-t_\tau x}{a t_\tau}\to \frac{y-tx}{at}.$$
Denote this limit by $z$. Then $z\in D_X$ (as $z_\tau\in D_X$ and $D_X$ is closed in the norm
topology and weak* closed in the dual case) and hence $y=t(x+az)\in D$. Moreover, $\|y\|\le 1$
(again by closedness of $D_X$), so $y\in D'$. This completes the proof.
\end{proof}

\begin{lemma}\label{balancedconstruction} Let $X$ be a Banach space,  $\varepsilon\in(0,1)$ and $(\xi_n)_{n\in\en}$  a sequence in $S_{X^*}$. Suppose that the following condition is satisfied:
\begin{itemize}
    \item [(i)] $\dist(\xi_n,\sp\{\xi_1,\dots,\xi_{n-1}\})>\varepsilon$ for all natural numbers $n\ge 2$.
\end{itemize}
Let $(a_n)_{n\in\en}$ be a sequence of real numbers satisfying the following conditions:
\begin{itemize}
\item[(ii)] $0<a_n<\frac{\varepsilon}{17}$ for each $n\in\en$;
\item[(iii)]  $a_n\to 0$;
\end{itemize}
Then the set
$$U:=\left(\frac12 B_{X^*}\cup\bigcup\{ t(\xi_n+a_n B_{X^*}) : t\in\cu\setminus\{0\} \}\right)\cap B_{X^*}$$ is a bounded
balanced norm-open set in $X^*$ satisfying  $\tilde U=U\subsetneqq B_{X^*}$.
   \end{lemma}

\begin{proof} It is clear that $U$ is an open balanced set and that $U \subset B_{X^*}$.
It remains to prove that $U\ne B_{X^*}$ and $\tilde U=U$.

For $n\in\en$ set
\begin{alignat*}{3}
D_n&:=\bigcup\{ t(\xi_n + a_n D_{X^*}) : t\in\cu \}\cap D_{X^*}, & D'_n&:=D_n\setminus\frac12 B_{X^*}, \\
C_n&:=\bigcup\{ t(\xi_n+a_n B_{X^*}) : t\in\cu\setminus\{0\}\}\cap B_{X^*}, & \qquad C'_n&:=C_n\setminus\frac12 B_{X^*}.
\end{alignat*}

The proof will continue by proving several consecutive claims.

\begin{itemize}
    \item[(a)] For each $n\in\en$ and $\eta\in D'_n$ there is some $t\in\cu$ with $|t|\in[\frac{1}{2(1+a_n)},\frac1{1-a_n}]$ with $\eta\in t(\xi_n+a_n D_{X^*})$.
\end{itemize}

Let $n\in\en$ and $\eta\in D'_n$ be arbitrary. By Lemma~\ref{cutcones} there is some $t\in\cu$ with $|t|\le\frac1{1-a_n}$ and $\eta\in t(\xi_n+a_n D_{X^*})$. We will check that $|t|$ satisfies also the lower bound:

Fix $\theta\in D_{X^*}$ such that $\eta=t(\xi_n+a_n\theta)$.
If $|t|<\frac{1}{2(1+a_n)}$, then
$$\|\eta\|\le |t|(\|\xi_n\|+a_n\|\theta\|)\le |t|(1+a_n)<\frac12,$$
so $\eta\notin D'_n$.

\begin{itemize}
    \item [(b)] $\dist(D'_n,D'_m)\ge \frac{\varepsilon}{4}$ for all $n, m \in \en$, $n\neq m$.
\end{itemize}

Let $m,n\in\en$ be such that $n\neq m$. Without loss of generality suppose $n>m$. Take $\eta_n\in
D'_n$ and $\eta_m\in D'_m$. By the already proved assertion (a) there are $\theta_m,\theta_n\in
D_{X^*}$ and  numbers $t_n,t_m$ with modulus in $[\frac{1}{2(1+a_n)},\frac1{1-a_n}]$ and
$[\frac{1}{2(1+a_m)},\frac1{1-a_m}]$ respectively such that
$$\eta_n=t_n(\xi_n+a_n\theta_n),\qquad \eta_m=t_m(\xi_m+a_n\theta_m).$$
Then
\begin{align*}
\|\eta_n-\eta_m\|&\ge |t_n|\cdot\left\|\xi_n-\frac{t_m}{t_n}\xi_m\right\|-|t_n|\cdot a_n \|\theta_n\|-|t_m|\cdot a_m\cdot\|\theta_m\|
\\ & \ge\frac1{2(1+a_n)}\cdot\varepsilon-
\frac{a_n}{1-a_n}-\frac{a_m}{1-a_m}
>\frac{\varepsilon}{2(1+\frac{\varepsilon}{17})} -2\frac{\frac{\varepsilon}{17}}{1-\frac{\varepsilon}{17}}
>\frac{\varepsilon}4,
\end{align*}
where we used the triangle inequality, assumptions (i) and (ii) and the fact that $\varepsilon<1$.

\begin{itemize}
        \item [(c)] $\wscl{C_n}=\ncl{C_n}=D_n$ and $\Int_{\|\cdot\|}(D_n)=C_n$ for all $n\in \en$.
   \end{itemize}

This follows from Lemma~\ref{cutcones}.

\begin{itemize}
\item[(d)] For each $n\in\en$ there is some $\zeta_n\in X^*$ with $\|\zeta_n\|=1$ such that $t\zeta_n\in D_n\setminus C_n$ for each $t\in\cu$ with
$0<|t|<1$.
\end{itemize}

Set $M=B_{X^*}\setminus\{0\}$. Then $C_n$ is an open subset of $M$ and, due to (c), $D_n\cap M$ is
the relative norm-closure of $C_n$ in $M$. As $D_n\cap M\subsetneqq M$ (this follows, for example,
from (b), as $D_m'\cap M$ is nonempty and disjoint from $D_n$ for each $m\in\en$, $m\ne n$) and
$M$ is obviously connected (note that $X^*$ is a nontrivial complex space, and so each two points
of $M$ can be joined in $M$ either by a segment or by an arc), we conclude that $(D_n\cap
M)\setminus C_n\ne\emptyset$. If we choose $\eta$ in this difference, it is enough to set
$\zeta_n=\frac{\eta}{\|\eta\|}$.

\begin{itemize}
\item[(e)] Let $n\in\en$ and $\eta\in C'_n$ be arbitrary. Then $\dist(\eta,X^*\setminus U)\le \frac{5a_n}{1-a_n}$.
\end{itemize}

Let $\zeta_n$ be the point given by (d). Then there is $\omega\in D_{X^*}$ and $s\in\cu$ with
$\zeta_n=s(\xi_n+a_n\omega)$. As $\|\zeta_n\|=1$, we get by the triangle inequality $|s|\in
[\frac1{1+a_n},\frac1{1-a_n}]$ (the computation is the same as that in Lemma~\ref{cutcones} and in
(a) above).

We further remark that $u\zeta_n\notin U$ for $|u|\ge\frac12$. Indeed, if $|u|\ge1$, then
$u\zeta_n\notin B_{X^*}$. Moreover, if $|u|\in[\frac12,1)$, then $\|u\zeta_n\|\ge\frac12$ and
$u\zeta_n\in D_n\setminus C_n$ by (d). Thus $u\zeta_n\notin D_m$ for $m\ne n$ (by (b)) and so
$u\zeta_n\notin U$.

Finally, choose any $\eta\in C'_n$. By (a) there is $\theta\in D_{X^*}$ and $t\in\cu$ with $|t|\in[\frac{1}{2(1+a_n)},\frac{1}{1-a_n}]$ such that $\eta=t(\xi_n+a_n\theta)$.

We will find $u\in\cu$ with $|u|\in[\frac12,1)$ such that
$|t-us|\le \frac{3a_n}{1-a_n}$. First let us find $\alpha\in\cu$ with $|\alpha|=1$ such that
$\alpha\frac st =|\frac st|$. Next we consider the following cases:

\begin{itemize}
  \item $|t|<\frac12|s|$. Set $u=\frac12\alpha$. Then $u$ has the required form and
  $$|t-us|=\left|t-\frac12\alpha s\right|=\frac12|s|-|t|\le \frac1{2(1-a_n)}-\frac1{2(1+a_n)}=\frac {a_n}{1-a_n^2}\le \frac{a_n}{1-a_n}.$$
    \item $\frac12|s|\le |t|<|s|$. Then we set $u=\frac ts$. This $u$ has the required form and $|t-us|=0$.
    \item $|t|>|s|$. Set $u=\frac{1-2a_n}{1-a_n}\alpha$. Then $u$ has the required form (it follows from (ii)) and
    $$|t-us|=\left|t- \frac{1-2a_n}{1-a_n}\alpha s\right|=|t|-\frac{1-2a_n}{1-a_n}|s|
    \le \frac{1}{1-a_n}-\frac{1-2a_n}{(1-a_n)(1+a_n)}=\frac{3a_n}{1-a_n^2}\le
    \frac{3a_n}{1-a_n}.$$
\end{itemize}

To conclude the proof of (e) observe that $u\zeta_n\notin U$ and
$$\|\eta-u\zeta_n\|=\|t(\xi_n+a_n\theta)-us(\xi_n+a_n\omega)\|
\le|t-us|+a_n(|t|+|us|)\le \frac{5a_n}{1-a_n}.$$

\begin{itemize}
    \item[(f)]  If $B\subset U$ is a $U$-bounded set, then there exists $n_0\in \en$ such that
$$B\subset \frac12 B_{X^*} \cup \bigcup_{n=1}^{n_0} C_n .$$
\end{itemize}

This follows immeadiately from (e) and the definition of a $U$-bounded set, as $\frac{5a_n}{1-a_n}\to0$.

\begin{itemize}
    \item[(g)]  $\Int_{\|\cdot\|}\left(\frac12D_{X^*}\cup\bigcup_{n\in\en} D_n\right)=U.$
\end{itemize}

The inclusion $\supset$ is obvious, let us prove the second one.
Let $\eta$ be in the set on the left-hand side.

If $\|\eta\|>\frac12$ then there exists $n\in\en$ such that $\eta\in D_n$. Choose  $s>0$ such that
$s<\min(\frac\varepsilon8,\|\eta\|-\frac12)$ and $\eta+s B_{X^*}$ is contained in the set on the
left-hand side. As $\eta+s B_{X^*}\cap D_m=\emptyset$ for all $m\neq n$ by (b), we get $\eta+
s_1B_{X^*}\subset D'_n$ hence $\eta\in C'_n$ by (c). Thus $\eta\in U$.

If $\|\eta\|=\frac12$, choose $s>0$ such that $\eta+s B_{X^*}$ is contained in the set on the left-hand side. Then $\eta'=(1+s)\eta$ also belongs to the set on the left-hand side and $\|\eta'\|>\frac12$. So, $\eta'\in U$ by the previous paragraph. As $U$ is balanced, we conclude that $\eta\in U$.

Finally, if $\|\eta\|<\frac12$, obviously $\eta\in U$.

\begin{itemize}
    \item[(h)] Conclusion: $\tilde U=U\subsetneqq B_{X^*}$.
\end{itemize}

That $U\subsetneqq B_{X^*}$ follows, for example, from (d) and (b), as $\frac34\zeta_n\in B_{X^*}\setminus U$ for each $n\in\en$.

If $B\subset U$ is U-bounded, let $n_0\in\en$ be the number provided by (f). Then
$$\wscl B\subset \wscl{\frac12 B_{X^*} \cup \bigcup_{n=1}^{n_0} C_n}
=\frac12 D_{X^*} \cup \bigcup_{n=1}^{n_0} D_n
\subset\frac12 D_{X^*} \cup \bigcup_{n\in\en} D_n.$$
Therefore
$$\tilde U\subset\frac12 D_{X^*} \cup \bigcup_{n\in\en} D_n.$$
As $\tilde U$ is norm-open, we conclude that
$$U\subset \tilde U\subset
\Int_{\|\cdot\|}(\frac12 D_{X^*} \cup \bigcup_{n\in\en} D_n)= U$$
by (g). This completes the proof.
\end{proof}

A sequence $(\xi_n)$ satisfying the assumptions of Lemma~\ref{balancedconstruction} can be constructed in the dual to any infinite-dimensional space by an easy application of the Riesz lemma.
The following lemma shows us that in some spaces we can choose such a sequence which is moreover weak* dense in $D_{X^*}$. In this way we prove Theorem~\ref{bal-cpcp}.

\begin{lemma}\label{og} Let $X$ be a separable Banach space and $\varepsilon'>2\varepsilon>0$ be
such that each nonempty relatively weak* open subset of $D_{X^*}$ has diameter greater than
$\varepsilon'$. Then there is a sequence $(\xi_n)$, $n\in\en$ in $X^*$ satisfying the following
conditions:
\begin{itemize}
    \item [(i)] $\|\xi_n\|=1$ for $n\in\en$.
    \item [(ii)] $\{\xi_n:n\in\en\}$ is weak* dense in $D_{X^*}$.
    \item [(iii)] $\dist(\xi_n,\sp\{\xi_1,\dots,\xi_{n-1}\})>\varepsilon$ for each $n\ge2$.
\end{itemize}
\end{lemma}

\begin{proof} As $X$ is separable, $(D_{X^*},w^*)$ is a metrizable
compact, and so it has a countable basis. Fix $U_n, n\in\en$ such a
basis consisting of nonempty sets.

Note that $X$ is necessarily infinite-dimensional, hence the sphere $S_{X^*}$ is weak* dense in $B_{X^*}$.

Choose $\xi_1\in U_1\cap S_{X^*}$ arbitrary.  Suppose we have constructed
$\xi_1,\dots,\xi_{n-1}$. We will choose $\xi_n\in U_n\cap S_{X^*}$ such that
$\dist(\xi_n,\sp\{\xi_1,\dots,\xi_{n-1}\})>\varepsilon.$

If we do that, then $\xi_n$'s clearly satisfy (i)--(iii). It remains
to show that the choice of $\xi_n$ is possible.

Set
$$F=(1+2\varepsilon)D_{X^*}\cap\sp\{\xi_1,\dots,\xi_{n-1}\}.$$
As $F$ is norm-compact, there are finitely many points $\eta_1,\dots,
\eta_m\in F$ such that
$$F\subset\{\eta_1,\dots,\eta_m\}+\left(\frac{\varepsilon'}2-\varepsilon\right) B_{X^*}.$$
Then we get
\begin{equation}\label{ssiitt}
F+\varepsilon D_{X^*}\subset \{\eta_1,\dots,\eta_m\}+\frac{\varepsilon'}2 D_{X^*}
\end{equation}
As each nonempty weak* open subset of $B_{X^*}$ has
diameter greater than $\varepsilon'$ and $U_n\ne\emptyset$, we get
$U_n\setminus \eta_1+\frac{\varepsilon'}2 D_{X^*}\ne\emptyset$. Repeating
the same argument $(n-2)$ times we get
$$U_n\setminus\{\eta_1,\dots,\eta_m\}+\frac{\varepsilon'}2 D_{X^*}\ne\emptyset.$$
Moreover, this set is clearly weak* open in $D_{X^*}$, therefore
we can choose $\xi_n$ in that set satisfying $\|\xi_n\|=1$. By
(\ref{ssiitt})
we get $\dist(\xi_n,F)>\varepsilon$. As for $\eta\in\sp\{\xi_1,\dots,\xi_n\}\setminus F$ we have $\|\eta\|>1+2\varepsilon$ and so $\|\xi_n-\eta\|>2\varepsilon$, we conclude that (iii) is
satisfied and the induction step is completed.
\end{proof}

Theorem~\ref{bal-cpcp} now follows immediately from Lemma~\ref{og} and Lemma~\ref{balancedconstruction}.

\section{Envelopes of unbounded convex sets}\label{s-convex}

In this section we study the behaviour of $\tilde U$ for unbounded open convex sets $U$. The main focus is on absolutely convex sets as, due to Proposition~\ref{espectrow}, it is natural to consider balanced sets. Anyway, most results of this section are of geometrical or topological nature and hold also for convex sets. Therefore we formulate those results in this more general setting which may be also interesting in itself.

We focus on the question when $\tilde{\tilde U}=\tilde U$. It is easy to see that the latter
equality holds whenever $U$ is boundedly regular. The converse does not hold as we will show in
the next section. We characterize those open convex sets which satisfy $\tilde{\tilde U}=\tilde U$
and describe a method of constructing counterexamples. In the next section we will use it to
provide a variety of examples of unbounded open absolutely convex sets satisfying $\tilde{\tilde
U}\ne\tilde U$.  This will be also applied to distinguishing the respective algebras of
holomorphic functions.

As we will deal with iteration of the $\widetilde\ $-envelope, we adopt the following notation. If $U\subset X^*$ is an open set and $\alpha$ is an ordinal, we define inductively
$$\Dt{\alpha}{U}=\begin{cases} U & \mbox{ if }\alpha=0, \\
\widetilde{\Dt{\beta}{U}} & \mbox{ if }\alpha=\beta+1, \\
\bigcup_{\beta<\alpha} \Dt{\beta}{U} & \mbox{ if $\alpha$ is limit.}
\end{cases}$$

It turns out that the $\widetilde\ $-envelope is closely related to the following operation: Let $A\subset X^*$ be any subset. Set
$$\Dd{1}{A}=\bigcup_{n\in\en}\wscl{A\cap n B_{X^*}},$$
i.e., $\Dd{1}{A}$ is the set of all limits of  weak* convergent bounded nets in $A$. If $X$ is
separable, it is just the set of all limits of weak* convergent sequences from $A$. We define
inductively $\Dd{\alpha}{A}$ for any ordinal $\alpha$ by the obvious way. We recall that, provided
$A$ is convex, by the Banach-Dieudonn\'e theorem $\Dd{1}{A}=A$ if and only if $A$ is weak* closed.

The relationship of the two operations is witnessed by the following theorem.

\begin{thm}\label{convextilde} Let $X$ be a Banach space and $U\subset X^*$ an open convex set. Then we have the following:
\begin{itemize}
    \item[(1)] $\tilde U=\Int_{\|\cdot\|} \Dd{1}{U}$.
    \item[(2)] $\Dt{\alpha}{U}=\Int_{\|\cdot\|} \Dd{\alpha}{U}$ for each ordinal $\alpha$.
    \item[(3)] Consider the following assertions.
    \begin{itemize}
    \item [(i)] $U$ is boundedly regular.
    \item[(ii)] $\tilde{\tilde U}=\tilde U$.
    \item[(iii)] $\tilde U=\icu$.
  \item[(iv)] $\wscl{U}=\ncl{\Dd{1}{U}}$.
\end{itemize}
Then the assertion (ii)--(iv) are equivalent and are implied by the assertion (i).
\end{itemize}
\end{thm}

Let us remark that the implication (ii)$\Rightarrow$(i) in the part (3) of the above theorem does not hold. An example is given in the next section.  Anyway,
already this theorem will enable us to construct absolutely convex open sets which are not boundedly regular by violating conditions (ii)--(iv).

To prove the theorem we need two lemmas.

\begin{lemma}\label{ffconvex} Let $X$ be a real or complex Banach space and $U$ an open convex set containg $0$. Then the sets
$$H_n=n B_{X}\cap \left(1-\frac1n\right) U, \qquad n\in\en,$$
form a fundamental family of convex $U$-bounded sets.
\end{lemma}

\begin{proof}
We will verify the three conditions from the respective definition.

First, we have to show that each $H_n$ is $U$-bounded. It is obviously bounded. Further, choose
$c>0$ with $cB_{X}\subset U$. Then
$$H_n+\frac cn B_{X}\subset \left(1-\frac1n\right) U+\frac1n U= U$$
as $U$ is convex. Thus $\dist(H_n, X\setminus U)\ge\frac cn$. This completes the proof that $H_n$
is $U$-bounded.

Further, let $A\subset U$ be $U$-bounded. We will find some $n\in\en$ with $A\subset H_n$. As $A$
is $U$-bounded, there is $M>0$ with $A\subset M B_{X}$ and $c>0$ with $A+c B_{X}\subset U$. Then
$$\left(1+\frac cM\right)A\subset A +\frac cM A\subset A + c B_{X}\subset U,$$
thus $A\subset \frac{1}{1+\frac cM} U$. It remains to take some $n\in\en$ such that $n\ge M$ and $\frac{1}{1+\frac cM}<1-\frac1n$. Then clearly $A\subset H_n$.

Finally, obviously
$$H_n+\left(\frac1{n}-\frac1{n+1}\right)B_{X}\subset H_{n+1}.$$
This completes the proof.
\end{proof}

\begin{lemma}\label{convexint}
Let $Y$ be a normed space and $U\subset Y$ a convex set such that $0\in\Int_{\|.\|} U$. Then
$$\Int_{\|.\|} U=\Int_{\|.\|} \overline{U}=\bigcup_{t\in[0,1)} tU.$$
\end{lemma}

\begin{proof} Let $p$ be the Minkonwski functional of $U$. As $0\in\Int_{\|.\|} U$, $p$ is continuous
and hence $A=\{x\in Y: p(x)<1\}$ is open and $B=\{x\in Y: p(x)\le 1\}$ is closed. Moreover, it is
clear that $A\subset U\subset B$, $A$ is dense in $B$ and $\Int_{\|.\|} B=A$. The announced
equalities now follow.
\end{proof}

{\bf Proof of Theorem~\ref{convextilde}.} Let $U\subset X^*$ be an open convex set. Without loss of generality we can suppose that $0\in U$.

(1) We have the following equalities:
\begin{align*}
\tilde U&=\bigcup_{n\in\en} \wscl{n B_{X^*}\cap \left(1-\frac1n\right) U}
=\bigcup_{n\in\en}\bigcup_{m\in\en} \wscl{m B_{X^*}\cap \left(1-\frac1n\right) U}
\\ &=\bigcup_{n\in\en} \Dd{1}{\left(\left(1-\frac1n\right) U\right)}
=\bigcup_{n\in\en} \left(1-\frac1n\right)\Dd{1}{ U}
=\Int_{\|\cdot\|}\Dd{1}{U}.
\end{align*}
Indeed, the first equality follows from Lemma~\ref{ffconvex} and the definition of $\tilde U$. The inclusion $\subset$ of the second equality is obvious, the converse one follows from the observation that
$m B_{X^*}\cap (1-\frac1n) U\subset k B_{X^*}\cap (1-\frac1k) U$ whenever $k\ge\max\{m,n\}$. The third equality follows immediately from the definitions,
the fourth one is obvious. The last one follows from Lemma~\ref{convexint} applied to the convex set $\Dd{1}{U}$.

Thus the assertion (1) is proved.

(2) If $\alpha=0$, the assertion is obvious, the case $\alpha=1$ is covered by the assertion (1). Suppose that $\alpha>1$ and that the assertion is true for all $\beta<\alpha$.

First suppose that $\alpha$ is isolated, i.e., $\alpha=\beta+1$ for some $\beta$.
Then
$$\Dt{\alpha}{U}=\widetilde{\Dt{\beta}{U}}
=\widetilde{\Int_{\|\cdot\|}\Dd{\beta}{U}}
=\Int_{\|\cdot\|}\Dd{1}{(\Int_{\|\cdot\|}\Dd{\beta}{U})}
\subset\Int_{\|\cdot\|}\Dd{1}{(\Dd{\beta}{U})}
=\Int_{\|\cdot\|}\Dd{\alpha}{U}.$$
Indeed, the first equality follows from the definitions, the second one from the induction hypothesis, the third one from the assertion (1). The rest is obvious.

In this way we have proved the inclusion $\subset$. To prove the converse one,
choose any $\xi\in \Int_{\|\cdot\|}\Dd{\alpha}{U}$. By Lemma~\ref{convexint} there is some $r\in[0,1)$ and $\eta\in \Dd{\alpha}{U}$ with $\xi=r\eta$.
Then there is $M>0$ such that $\eta\in\wscl{\Dd{\beta}{U}\cap M B_{X^*}}$.
Fix some $s\in(r,1)$. Then we have
$$s\eta\in\wscl{s\Dd{\beta}{U}\cap sM B_{X^*}}\subset \widetilde{\Dd{\beta}{U}}.$$
Indeed, this follows from Lemma~\ref{ffconvex} if we choose $n\in\en$ with $n>sM$ and $1-\frac1n>s$.
So,
$$\xi=r\eta=\frac rs \cdot s\eta\in\frac rs\widetilde{(\Dd{\beta}{U})}
=\widetilde{\left(\frac rs\Dd{\beta}{U}\right)}
\subset \widetilde{(\Int_{\|\cdot\|}\Dd{\beta}{U})}
=\widetilde{(\Dt{\beta}{U})}=\Dt{\alpha}{U}.$$
The first two equalities are trivial. The following relation is proved above. The third equality
is trivial, the next inclusion follows from Lemma~\ref{convexint}, the next equality follows from
the induction hypothesis. The last one is just the definition.

Next suppose that $\alpha$ is limit. Then
\begin{align*}
\Dt{\alpha}{U}&=\bigcup_{\beta<\alpha} \Dt{\beta}{U}
=\bigcup_{\beta<\alpha}\Int_{\|\cdot\|} \Dd{\beta}{U}
=\bigcup_{\beta<\alpha}\bigcup_{t\in[0,1)} t\Dd{\beta}{U}
\\ &=\bigcup_{t\in[0,1)} t\left(\bigcup_{\beta<\alpha}\Dd{\beta}{U}\right)
=\bigcup_{t\in[0,1)} t\Dd{\alpha}{U}
=\Int_{\|\cdot\|} \Dd{\alpha}{U},
\end{align*}
where we used the definitions, induction hypothesis and Lemma~\ref{convexint}.

This completes the proof of (2).

(3) The implication (i)$\Rightarrow$(ii) is obvious.

(ii)$\Rightarrow$(iii) Suppose that $\tilde{\tilde U}=\tilde U$. Then for each ordinal $\alpha\ge
1$ we have $\tilde U=\Dt{\alpha}{U}$. It follows from the Banach-Dieudonn\'e theorem that there is
an ordinal $\alpha$ such that $\Dd{\alpha}{U}=\wscl{U}$. Then we have by (2) that
$$\tilde U=\Dt{\alpha}{U}=\Int_{\|\cdot\|}\Dd{\alpha}{U}=\icu,$$
which completes the proof.

(iii)$\Rightarrow$(iv) Suppose (iii) holds. By Lemma~\ref{convexint} we get that
$$\tilde U=\bigcup_{t\in[0,1)}t\wscl{U},$$
hence $\tilde U$ is norm-dense in $\wscl{U}$. As $\tilde U\subset\Dd{1}{U}\subset\wscl{U}$ by (1), we conclude that $\Dd{1}{U}$ is a norm-dense subset of $\wscl{U}$, thus (iv) is valid.

(iv)$\Rightarrow$(iii) Suppose (iv) holds.
Then
$$\icu=\Int_{\|\cdot\|}\ncl{\Dd{1}{U}}=\Int_{\|\cdot\|}\Dd{1}{U}=\tilde U.$$
Indeed, the first equality follows from the assumption (iv), the second one from Lemma~\ref{convexint} and the last one from the assertion (1).

(iii)$\Rightarrow$(ii) is obvious due to (2) as
$$\tilde U\subset\tilde{\tilde U}=\Int_{\|\cdot\|}\Dd{2}{U}\subset\icu.$$

\hfill\qed\medskip

Now we turn to the method of constructing counterexamples. The key tool is the following proposition.

\begin{prop}\label{Aplusball} Let $X$ be a Banach space, $A\subset X^*$ a convex set and $r>0$. Then $U=A+r B_{X^*}$ is an open convex set and for each ordinal $\alpha$ we have
$$\Dt{\alpha}{U}=\Dd{\alpha}{A}+r B_{X^*}.$$
\end{prop}

To prove this proposition we need the following easy lemma.

\begin{lemma}\label{distance} Let $Y$ be a normed space, $A\subset Y$ a nonempty convex set, $x\in Y$ and $d>0$ such that $\dist(x,A)\ge d$. Then for each $\delta\in(0,1)$  there is $y\in S_Y$ such that for each $t>0$ we have
$\dist(x+ty)\ge d+(1-\delta)t$.
\end{lemma}

\begin{proof} Suppose $Y$, $A$, $x$ and $d$ are as in the assumptions. Set
$B=x+ d B_Y$. Then $B$ is an open convex set and $A\cap B=\emptyset$. Hence by the Hahn-Banach separation theorem there is some $\xi\in S_{Y^*}$ and $c\in\er$ such that
$$\sup \Re\xi(A)\le c\le\inf\Re\xi(B).$$
Note that $\inf\Re\xi(B)=\Re\xi(x)-d$, so $\Re\xi(x)\ge c+d$.

If $\delta\in(0,1)$ is given, find $y\in S_Y$ such that $\xi(y)\in\er$ and $\xi(y)>1-\delta$. We claim this $y$ has the required properties. Indeed, let $t>0$ and $z\in A$ arbitrary. Then
\begin{align*}
\|(x+ty)-z\|&\ge |\xi(x+ty-z)|\ge\Re\xi(x+ty-z)\\&=
\Re\xi(x)+t\xi(y)-\Re\xi(z)\ge (c+d)+t(1-\delta)-c
=d+t(1-\delta).
\end{align*}
This finishes the proof.
\end{proof}

{\bf Proof of Proposition~\ref{Aplusball}.} It is enough to prove the proposition for $\alpha=1$.
The general case is then immediate by transfinite induction. We will prove that
$$\tilde U=\Dd{1}{A}+r B_{X^*}.$$

That $\Dd{1}{A}+r B_{X^*} \subset \tilde U$ is obvious, since
$$A\cap n B_{X^*}+sB_{X^*}+(r-s)B_{X^*}=A\cap n B_{X^*}+rB_{X^*}\subset U,$$
for all $n\in\en$ and $0<s<r$, so $A\cap n B_{X^*}+sB_{X^*}$ is $U$-bounded.
Finally, it is clear that weak* closure of these sets cover $\Dd{1}{A}+r B_{X^*}$.

Conversely, suppose $\xi\in\tilde U$. Then there is $M\subset\tilde U$ which is $U$-bounded and
$\xi\in\wscl M$. First let us prove that there is some $s\in(0,r)$ with $M\subset A+s B_{X^*}$.

As the distance of $M$ to the complement of $U$ is positive, there is some $c>0$ with $M+c B_{X^*}\subset U$. Choose any $\eta\in M$ and set $d=\dist(\eta,A)$.
Suppose that $d>0$. By Lemma~\ref{distance} there is $\zeta\in S_{X^*}$
such that for each $t>0$ we have $\dist(\eta+t\zeta,A)\ge d+\frac t2$.
As $\eta+ c B_{X^*}\subset U$, we get
$\dist(\eta+t\zeta,A)< r$ whenever $t\in(0,c)$. Therefore
$r\ge d+\frac c2$, so $d\le r-\frac c2$. It follows that $r-\frac c2\ge 0$ and that $M\subset A+(r-\frac c2) D_{X^*}$, so $M\subset A+s B_{X^*}$ for any $s\in(r-\frac c2,r)$.

So fix some $s\in(0,r)$ with $M\subset A+s B_{X^*}$. The set $M$ is also bounded, so there is some $R>0$ with $M\subset R B_{X^*}$. Then clearly
$$M\subset (A\cap (R+s) B_{X^*}) + s B_{X^*},$$
hence
$$\xi\in\wscl{M}\subset \wscl{(A\cap (R+s) B_{X^*}) + s B_{X^*}}
\subset \wscl{(A\cap (R+s) B_{X^*})} + s D_{X^*}
\subset\Dd{1}{A}+r B_{X^*}.$$
This completes the proof.
\qed\medskip

We finish this section by the following proposition which shows how to construct examples of (absolutely) convex sets which are not boundedly regular.

\begin{prop} Let $X$ be a Banach space. The following assertions are equivalent:
\begin{itemize}
    \item[(i)] For each (absolutely) convex open set $U\subset X^*$ we have $\tilde{\tilde U}=\tilde U$.
    \item[(ii)] For each (absolutely) convex set $A\subset X^*$ we have $\wscl{A}=\ncl{\Dd{1}{A}}$.
\end{itemize}
\end{prop}

\begin{proof}
The implication (ii)$\Rightarrow$(i) follows immediately from Theorem~\ref{convextilde}.(3).

Let us show (i)$\Rightarrow$(ii). Suppose (ii) does not hold. So there is an (absolutely) convex
set $A\subset X^*$ with  $\wscl{A}\supsetneqq\ncl{\Dd{1}{A}}$. So, we can fix
$\xi\in\wscl{A}\setminus\ncl{\Dd{1}{A}}$. Further, find $r>0$ such that $\dist(\xi,\Dd{1}{A})>r$
and set $U=A+r B_{X^*}$. Then $U$ is (absolutely) convex and open and, moreover, $\tilde
U=\Dd{1}{A}+r B_{X^*}$ by Proposition~\ref{Aplusball}, so $\xi\notin\tilde U$. On the other hand,
$\xi\in\wscl{A}$, so there is an ordinal $\alpha$ with $\xi\in \Dd{\alpha}{A}$ and therefore also
$\xi\in\Dt{\alpha}{U}$. Thus, $\Dt{\alpha}{U}\ne\tilde U$, hence $\tilde{\tilde U}\ne\tilde U$.
\end{proof}

\section{Absolutely convex sets which are not boundedly regular}\label{s-examples}

In this section we collect examples of absolutely convex sets $U$ which are not boundedly regular.
Most of them even do not satisfy $\tilde{\tilde U}=\tilde U$, but we also give  one example
satisfying $\tilde{\tilde U}=\tilde U$.

The method of construction consists in using Proposition~\ref{Aplusball} with $A$ being a linear subspace of $X^*$. Recall that a Banach space $X$ is called quasireflexive if the quotient $X^{**}/X$ has finite dimension. We have the following well-known dichotomy.

\begin{prop}\label{quasir-subs}
\begin{itemize}
    \item[(1)] Let $X$ be a Banach space. Then $X$ is quasireflexive  if and only if $\Dd{1}{E}=\wscl{E}$ for each linear subspace $E\subset X^*$.
    \item[(2)] Let $E$ be a separable space which is not quasireflexive. Then for each countable ordinal $\alpha$ there is a linear subspace $E\subset X^*$ such that $\Dd{\alpha+1}{E}=X^*$ and $\Dd{\alpha}{E}$ is contained in a proper closed hyperplane in $X^*$.
\end{itemize}
\end{prop}

\begin{proof}
(1) This follows from \cite{Pet} (the only if part) and \cite{DL} (the if part).

(2) Let $X$ be a separable Banach space which is not quasireflexive and $\alpha$ be a countable ordinal. By \cite[Theorem]{Ostro1} there is a linear subspace $E\subset X^*$ such that $\Dd{\alpha}{E}\subsetneqq\Dd{\alpha+1}{E}=X^*$. Moreover, it follows from the proof that $\Dd{\alpha}{E}$ is contained in a proper closed hyperplane in $X^*$. Indeed, one can take $E$ to be of the form
$K(g_0,\alpha,A)$ in the notation of \cite[Lemma 2]{Ostro1}. By the assertion 1) of the quoted lemma we have $\Dd{\alpha}{E}\subset \Ker g_0$. As $g_0\in X^{**}$, the proof is completed.
\end{proof}

As a consequence we get the following theorem.

\begin{thm}\label{t-qr-s} Let $X$ be a  separable Banach space which is
not quasireflexive. Then,
for each  $\alpha<\omega_1$ there is an (unbounded) absolutely convex open set $U$ in $X^\ast$ such that
$\Dt{\beta}{U}\varsubsetneq \Dt{\beta+1}{U}$ for every ordinal $\beta\le\alpha$ and $\Dt{\alpha+1}{U}=X^*$.
\end{thm}

\begin{proof}
Let $X$ be a separable Banach space which is not quasireflexive and $\alpha$ be a countable ordinal. Let $E\subset X^*$ be a linear subspace given by Theorem~\ref{quasir-subs} and $r>0$ be arbitrary. Set $U=E+rB_{X^*}$.
Then $U$ is an open absolutely convex set. By Proposition~\ref{Aplusball} we have
$\Dt{\beta}{U}=\Dd{\beta}{E}+r B_{X^*}$ for each ordinal $\beta$.
In particular,
$$\Dt{\alpha+1}{U}=\Dd{\alpha+1}{E}+rB_{X^*}=X^*$$
and
$$\Dt{\alpha}{U}=\Dd{\alpha}{E}+rB_{X^*}\subset H+rB_{X^*}\subsetneqq X^*,$$
where $H$ is a closed hyperplane in $X^*$.
This completes the proof.
\end{proof}

Even if the above theorem yields plenty of counter-examples, we are going to give an explicite one.

\begin{example2} Let
$$A=\{x\in\ell_\infty(\en\times\en): \forall k\in\en: \{ n\in \en:
x(k,n)\ne k x(k,1)\} \mbox{ is finite}\}.$$
Then $A$ is a linear subspace of $\ell_\infty(\en\times\en)=\ell_1(\en\times\en)^*$.
Then
$$\{x\in\ell_\infty(\en\times\en):\{(k,n): x(k,n)\ne0\} \mbox{ is
finite}\}\subset \Dd{1}{A} \subset\{x\in\ell_\infty(\en\times\en):
\lim_{k\to\infty}x(k,1)=0\}$$

Indeed, let $x$ be in the set on the leftt-hand side. Then there is
some $N\in\en$ such that $x(k,n)=0$ whenever $k>N$ or $n>N$. Define
a sequence $x_m$ as follows:
$$x_m(k,n)=\begin{cases} x(k,n), & k\le N,n\le N, \\ k x(k,1),& k\le
N, n\ge N+m, \\ 0, & \mbox{otherwise}.\end{cases}$$
Then $x_m$ is a bounded sequence of elements of $A$ which weak*
converges to $x$. This proves the first inclusion.

Let us show the second one. Let $x\in \Dd{1}{A}$. Then there is $N\in\en$ such that $x\in\wscl{\{y\in A:\|y\|\le N\}}$. Let $\varepsilon>0$ and $k\in\en$ such that $|x(k,1)|>\varepsilon$. Then there is $y\in A$ with $\|y\|\le N$ such that
$|y(k,1)|>\varepsilon$. By the definition of $A$ there is some (in fact many) $n\in\en$ with $y(k,n)=k y(k,1)$, so $|y(k,n)|>k\varepsilon$. Thus $k\varepsilon<N$, so $k<\frac N\varepsilon$.

It follows that $|x(k,1)|\le\varepsilon$ if $k\ge\frac N\varepsilon$. Therefore $\lim x(k,1)=0$, which completes the proof of the second inclusion.

From the first inclusion we get $\Dd{2}{A}=\ell_\infty(\en\times\en)$, the second one implies that $\Dd{1}{A}$ is contained in a closed hyperplane.
Thus if we set, for example, $U=A+B_{\ell_\infty(\en\times\en)}$ we obtain an open absolutely convex set with $\tilde U\subsetneqq\tilde{\tilde U}$.
\end{example2}

The following example is a modification of the previous one. It will enable us to distinguish boundedly regular sets and sets satisfying $\tilde{\tilde U}=\tilde U$.

\begin{example}\label{exa-ell1} There is a linear subspace $A\subset \ell_1=c_0^*$ such that $\Dd{1}{A}$ is a proper norm-dense subset of $\ell_1$.
\end{example}

\begin{proof} We will consider $\ell_1(\en\times\en)$ and set
$$A=\left\{x\in\ell_1(\en\times\en) : (\forall k\in\en)\left(x(k,1)=\frac1k\sum_{n=2}^\infty x(k,n)\right) \right\}.$$
Then $A$ is clearly a linear subspace of $\ell_1(\en\times\en)$. First we will show that $\Dd{1}{A}$ is norm-dense in $\ell_1(\en\times\en)$. To this end it is enough to show that each finitely supported vector belongs to $\Dd{1}{A}$.

Let $x\in\ell_1(\en\times\en)$ be finitely supported. Then there is
some $N\in\en$ such that $x(k,n)=0$ whenever $k>N$ or $n>N$. Define
a sequence $x_m$ as follows:
$$x_m(k,n)=\begin{cases} kx(k,1)-\sum_{j=2}^N x(k,j), & k\le N,n=N+m, \\  x(k,n)&  \mbox{otherwise}.\end{cases}$$
Then clearly $x_m\in A$ and the sequence $x_n$ converges to $x$ pointwise on $\en\times\en$. Moreover,
$$\|x_m\|\leq \|x\|+\sum_{k=1}^N\left|kx(x,1)-\sum_{j=2}^N x(k,j)\right|\le (1+N)\|x\|,$$
so the sequence $x_m$ is bounded. Hence the sequence $x_m$ weak* converges to $x$ and so $x\in\Dd{1}{A}$.

Finally, we will show that $\Dd{1}{A}$ is a proper subset of $\ell_1(\en\times\en)$. Namely, the element $x$ defined as
$$x(k,n)=\begin{cases} \frac1{k^2}, & n=1, \\ 0 & \mbox{otherwise}\end{cases}$$
belongs to $\ell_1(\en\times\en)\setminus \Dd{1}{A}$. Indeed, suppose that
$x\in\Dd{1}{A}$. Then there is $M>0$ such that $x\in\wscl{A\cap MB_{\ell_1(\en\times\en)}}$.
Choose $N\in\en$ such that $\sum_{k=1}^N\frac1{2k}>M$.
We can find $y\in A$ with $\|y\|<M$ such that $|y(k,1)|>\frac1{2k^2}$ for $k=1,\dots,N$. As $y\in A$, we get
$$\|y\|=\sum_{k,n\in\en}|y(k,n)|\ge\sum_{k=1}^N\left|\sum_{n=2}^\infty y(k,n)\right|
=\sum_{k=1}^N k|y(k,1)|>\sum_{k=1}^N\frac{1}{2k}>M,$$
a contradiction.
\end{proof}

The following proposition contains a characterization of boundedly regular sets
among the sets of the form considered in this section.

\begin{prop}\label{Apdd} Let $X$ be a Banach space and $A\subset X^*$ a weak*-dense linear subspace. Let $r>0$ be arbitrary and $U=A+rB_{X^*}$. Then $U$ is boundedly regular if and only if $\Dd{1}{A}=X^*$, i.e. $A$ is norming.
\end{prop}

\begin{proof} First suppose that $\Dd{1}{A}=X^*$.
Then $\wscl{A\cap B_{X^*}}$ contains $c B_{X^*}$ for some $c>0$. (This is an easy and well-known consequence of the Baire category theorem.)
Let $C\subset \tilde U$ be $\tilde U$-bounded. As $\tilde U=X^*$ (by Proposition~\ref{Aplusball}), it means just that $C$ is bounded, i.e., $C\subset MD_{X^*}$ for some $M>0$. As
$$M D_{X^*}=\frac Mc c D_{X^*}= \frac Mc \wscl{cB_{X^*}}
\subset\frac Mc \wscl{A\cap B_{X^*}}=\wscl{A\cap \frac Mc B_{X^*}}.$$
As $A\cap\frac Mc B_{X^*}$ is $U$-bounded, we conclude that $U$ is boundedly regular.

Conversely, suppose that $U$ is boundedly regular. Then necessarily $\tilde{\tilde U}=\tilde U$ and so $\tilde U= X^*$. (As $\wscl A=X^*$, there is an ordinal $\alpha$ such that $\Dd{\alpha}{A}=X^*$. Then $\Dt{\alpha}{U}=X^*$ as well by Proposition~\ref{Aplusball}.) The set $2r B_{X^*}$ is then $\tilde U$-bounded. As $U$ is boundedly regular, there is a $U$-bounded set $C\subset U$
such that $2r B_{X^*}\subset \wscl C$. $C$ is necessarily bounded, so there is some $M>0$ such that $C\subset (A\cap M B_{X^*})+r B_{X^*}$. So, we have
$$2r B_{X^*}\subset\wscl{(A\cap M B_{X^*})+r B_{X^*}}=
\wscl{A\cap M B_{X^*}} + r D_{X^*}.$$
For $x\in X$ set
$$p(x)=\sup\{|\xi(x)|: \xi\in \wscl{A\cap M B_{X^*}} \}.$$
Then $p$ is clearly a seminorm on $X$ such that $p(x)\le M\|x\|$ for each $x\in X$.

Fix $x\in X$ with $\|x\|=1$. There is $\xi\in S_{X^*}$ with $\xi(x)=1$.
Then $2r\xi= \eta + r\theta$ with $\eta\in\wscl{A\cap MB_{X^*}}$ and  $\theta\in D_{X^*}$. So,
$$p(x)\ge|\eta(x)|=|2r\xi(x)-r\theta(x)|\ge 2r |\xi(x)|-r|\theta(x)|\ge 2r-r=r.$$
It follows that $p(x)\ge r\|x\|$ for $x\in X$, so $p$ is an equivalent norm on $X$. By the bipolar theorem the respective dual unit ball is $\wscl{A\cap MB_{X^*}}$. Thus the latter set is a unit ball of an equivalent norm on $X^*$, so it contains $c B_{X^*}$ for some $c>0$. It follows that $\Dd{1}{A}=X^*$.
\end{proof}

Now we are ready to give the announced example:

\begin{example} There is an open absolutely convex set $U\subset \ell_1=c_0^*$ such that $\tilde{\tilde U}=\tilde U$ but $U$ is not boundedly regular.
\end{example}

\begin{proof} Let $A$ be the linear subspace of $\ell_1$ provided by Example~\ref{exa-ell1} and set $U=A+ B_{\ell_1}$. As $\Dd{1}{A}$ is norm dense in $\ell_1$, $\tilde U=\ell_1$ and so $\tilde{\tilde U}=\tilde U$.
Finaly, as $\Dd{1}{A}\ne\ell_1$, $U$ is not boundedly regular by Proposition~\ref{Apdd}.
\end{proof}

Finally, we use Proposition~\ref{Apdd} to partially extend Theorem~\ref{t-qr-s} to non-separable setting.

\begin{thm} Let $X$ be a (not necessarily separable) Banach space which is not quasireflexive. Then there is an open absolutely convex set $U\subset X^*$ which is not boundedly regular.
\end{thm}

\begin{proof} By \cite{DL} there is a weak* dense linear subspace $A\subset X^*$
such that $\Dd{1}{A}\ne X^*$. Then $U=A+ B_{X^*}$ is the required example due to Proposition~\ref{Apdd}.
\end{proof}

\section{Distinguishing algebras of holomorphic functions}\label{s-algebras}

The aim of this section is to prove the following theorem, which can be viewed as a kind of converse of Proposition~\ref{espectrow} for absolutely convex sets.

\begin{thm}\label{t-ext}
Let $X$ be a complex Banach space and $U\subset X^*$ an absolutely convex open set. Then each $f\in H_{w^\ast u}(U)$ can be extended to some $\tilde f\in H_{w^\ast u}(\tilde U)$ if and only if $U$ is boundedly regular.
\end{thm}

Let $(U_n)_{n\in\en}$ be a fundamental sequence of $U$-bounded sets consisting of absolutely convex sets (one can use Lemma~\ref{ffconvex}) and set
$\mathcal{B}=(\wscl{U_n})_{n\in\en}$.

The if part of the above theorem follows from Proposition~\ref{espectrow}, as for boundedly regular $U$ the algebras $H_{\mathcal{B}w^\ast u}(\tilde U)$ and $H_{w^\ast u}(\tilde U)$ coincide.

To prove the only if part we use the following proposition.

\begin{prop}  Let $X$ be a complex Banach space and  $U$ be an absolutely convex open set in $X^\ast$ which is not boundedly regular. Then $H_{w^\ast u}(\tilde U)$ is a proper subset of $H_{\mathcal{B}w^\ast u}(\tilde U)$ and  the topology of $H_{w^\ast u}(\tilde U)$ is strictly stronger than the topology inherited from $H_{\mathcal{B}w^\ast u}(\tilde U)$. \end{prop}

\begin{proof}
As $U$ is not boundedly regular, there exists a subset $E$ of $ \tilde U$ such that is  $\tilde U$-bounded but  is not contained in any $\wscl{U_n}$, for $n=1,2,...$. Hence we can take $x^\ast_n\in E\setminus \wscl{U_n}$. By applying the bipolar theorem, we can find an element $x_n\in X$ such that $ |x^\ast_n(x_n)|>1$ and  $|x^\ast(x_n)|\leq 1$ for all $x^\ast\in\wscl{U_n}$.

We choose a $p_n\in \na$ such that
\begin{equation*}|x^\ast_n(x_n)|^{p_n}>n,\end{equation*}
for each $n\in \na$. We take the $p_n$-homogeneous polynomial
$P_n:X^\ast\to \cn$ defined by
$$P_n(x^\ast)=(x^\ast(x_n))^{p_n}.$$
Then $P_n\in H_{w^\ast u}(\tilde U)$ and  $P_n\in H_{\mathcal{B}w^\ast u}(\tilde U)$ for all $n\in \na$.

 But $\sup\{|P_n(x^\ast)| : x^\ast \in E\}\geq n$ for all $n$ and $E$ is a  $ \tilde U$-bounded set. Hence the sequence $(P_n)$ is unbounded in $H_{w^\ast u}(\tilde U)$, but on the other hand  $\sup\{|P_n(x^\ast)| : x^\ast \in \wscl{U_p}\}\leq   \sup\{|P_n(x^\ast)| : x^\ast \in \wscl{U_n}\}\leq 1$ for all $n\geq p$ and any continuous polynomial is bounded on bounded subsets of $X^\ast$. Thus sequence $(P_n)$ is bounded in $H_{\mathcal{B}w^\ast u}(\tilde U)$.

Since the inclusion $i:H_{w^\ast u}(\tilde U)\longrightarrow H_{\mathcal{B}w^\ast u}(\tilde U)$ is always continuous, by the open mapping theorem, we have that  $H_{w^\ast u}(\tilde U)$ is a proper subspace of $H_{\mathcal{B}w^\ast u}(\tilde U)$ and the topology of $H_{w^\ast u}(\tilde U)$ is strictly stronger than the restriction of the topology $H_{\mathcal{B}w^\ast u}(\tilde U)$ to that subspace.
\end{proof}

Now we are ready to prove the only if part of Theorem~\ref{t-ext}.
If $U$ is not boundedly regular, by the previous proposition there is some
$f\in H_{\mathcal{B}w^\ast u}(\tilde U)\setminus H_{w^\ast u}(\tilde U)$. Let $g$ be the restriction of $f$ to $U$. Suppose that $\tilde g\in H_{w^\ast u}(\tilde U)$ is an extension of $g$. As we have $\tilde g\in H_{\mathcal{B}w^\ast u}(\tilde U)$ as well, $\tilde g=f$ by Proposition~\ref{espectrow} (by the uniqueness of the extension). But this is a contradiction as $f\notin H_{w^\ast u}(\tilde U)$.

\section{Final remarks and open problems}

In this section we collect some questions we do not know the answers. We begin by questions on balanced sets.

\begin{ques} Is there a bounded balanced open set in a dual Banach space which is not boundedly regular? In particular, is $\tilde{\tilde{U}}=\tilde{U}$ for any
bounded balanced open set? \end{ques}

In \cite{CaGaMa} an example of an unbounded balanced set which is not boundedly regular is given. Our example in Section~\ref{s-balanced} shows that balanced sets have very different behaviour from convex ones but still is boundedly regular.

\begin{ques} Is Theorem~\ref{t-ext} valid also for balanced sets?
\end{ques}

Note that in the proof of Theorem~\ref{t-ext} the convexity was essentially used due to the use of the bipolar theorem.

We continue by questions on convex sets:

\begin{ques} Let $X$ be a quasireflexive Banach space.
\begin{itemize}
    \item[(a)] Is $\wscl{A}=\Dd{1}{A}$ for each (absolutely) convex set $A\subset X^*$?
    \item[(b)] Is each (absolutely) convex open subset of $X^*$ boundedly regular?
\end{itemize}
\end{ques}

Note that the question (a) has positive answer if $A$ is a linear subspace. Anyway the respective proof strongly uses linearity and it seems not to be clear how to adapt it to the (absolutely) convex case.

We also remark that both questions have positive answer if $X$ is reflexive. Indeed, in this case weak* topology on $X^*$ conincides with the weak one.
Moreover, the weak closure of any convex set equals to its norm closure.
So, $\tilde U=U$ for each open convex set $U\subset X^*$.

\begin{ques} Let $X$ be a Banach space, $A\subset X^*$ an (absolutely) convex open set. Is it true that $A+r B_{X^*}$ is boundedly regular for each $r>0$ if and only if $\wscl{A}=\Dd{1}{A}$?
\end{ques}

Proposition~\ref{Apdd} shows that it is true if $A$ is a weak* dense linear subspace.

\begin{ques} For which Banach spaces $X$ there is a linear subspace $A\subset X^*$ such that $\Dd{1}{A}$ is a proper norm-dense subset of $X^*$? Is it true whenever $X$ is not quasireflexive?
\end{ques}

Note that we have just one example of such a subspace for $X=c_0$. It seems not to be clear how to adapt it even for $X=\ell_1$.

\section*{Acknowledgement}

We are grateful to Ji\v{r}\'{\i} Spurn\'y for helpful discussions on the subject.


\begin{thebibliography}{99}

\bibitem{AcosBe} M. D. Acosta and J. Becerra Guerrero,  Slices in the unit ball of the
symmetric tensor product of $C(K)$ and $L_1(\mu)$, {\it Ark. Mat.} {\bf 47} (2009), no. 1, 1--12.

\bibitem{AR}  R. M. Aron and P. Rueda, Homomorphisms on spaces of weakly
continuous holomorphic functions, {\it Arch. Math.} {\bf 73} (1999),
430--438.

\bibitem{Barroso} J.A. Barroso, {\it Introduction to holomorphy}, North-Holland
Mathematics Studies, 106. North-Holland Publishing Co., Amsterdam, 1985.

\bibitem{BLR} J. Becerra Guerrero, G. L\'opez P\'erez
and A. Rodr\'{\i}guez-Palacios, Relatively Weakly Open Sets in Closed Balls of C*-Algebras, {\it J. London Math. Soc.} {\bf 68} (2003), no. 3, 753--761.

\bibitem{CCG}  F. Cabello S\'{a}nchez, J. Castillo and R. Garc\'\i a,
Polynomials on dual-isomorphic spaces, {\it Ark. Mat.} {\bf 38}
(2000), 37--44.


\bibitem{CaGaMa} D. Carando, G. Garc\'{\i}a, M. Maestre,  Homomorphisms and composition operators
on algebras of analytic functions of bounded type, {\it Advances in
Mathematics} {\bf 197} (2005), 607-629.

\bibitem{DaGa} A. Davie and T. Gamelin,  A theorem on polynomial-star
approximation, {\it Proc. Amer. Math. Soc.} {\bf 106} (1989), 351-356.

\bibitem{DL}
{ W.J.Davis and J.Lindenstrauss},  On total nonnorming subspaces, {\it Proc. Amer. Math. Soc.},
{\bf 31} (1972), 109--111.

\bibitem{DGZ} { R.Deville, G.Godefroy and V.Zizler},
{\it Smoothness and renorming of Banach spaces}, Pitman Monographs, 1993.



\bibitem{DD} J.C. D\'{\i}az and S. Dineen, Polynomials on stable spaces, {\it Ark. Mat.} {\bf 36}
(1998), 87--96.

\bibitem{D}  S. Dineen, {\it Complex Analysis on Infinite Dimensional Spaces},
Spinger Monographs in Mathematics, Springer-Verlag, London, 1999.



\bibitem{LZ}  S. Lassalle and I. Zalduendo, To what extent does the
dual Banach space $E'$ determine the polynomials over $E$?, {\it
Ark. Mat.} {\bf 38} (2000), 343-354.


\bibitem{Ostro1}  M.I. Ostrovskii, $w^*$-derived sets of transfinite order of subspaces of dual Banach
spaces, Dokl. Akad. Nauk Ukrain. SSR, 1987, N.10, pp. 9--12 (in Russian and Ukrainian); (An English version of this paper can be found on the web at http://front.math.ucdavis.edu).



\bibitem{Pet}
{ Yu.I. Petunin}, Conjugate Banach spaces containing  subspaces  of zero  characteristic, {\it
Dokl. Akad.  Nauk  SSSR}, {\bf 154} (1964), 527--529 (in Russian); English transl. in: {\it Soviet
Math. Dokl.}, {\bf 5} (1964), 131--133.


\bibitem{phelps}
{ R.R. Phelps}, {\it  Convex Functions, Monotone Operators and Differentiability},
    Lecture Notes in Mathematics, Volume 1364/1989, 1993

\end{thebibliography}
 \end{document}